\documentclass[a4paper,11pt, twoside]{article}
\usepackage{a4}
\usepackage{geometry} 
\usepackage{amsmath}
\usepackage{amssymb}
\usepackage{amsthm}
\usepackage{graphicx}
\usepackage{caption,subcaption}
\usepackage{multirow}
\usepackage{xstring}
\usepackage[utf8]{inputenc}
\usepackage{amstext,amsbsy,amsopn,eucal,enumerate}
\usepackage{tikz}
\usepackage{pgfplots}
\usepackage{tabularx}

\def\R{\mathbb{R}}

\newcommand{\expn}{\operatorname{e}}

\newcommand{\diag}{\operatorname{diag}}

\newcommand{\beq}{\begin{equation}}
\newcommand{\eeq}{\end{equation}}
\newcommand {\mat}      [1] {\left[\begin{array}{#1}}
\newcommand {\rix}          {\end{array}\right]}
\newcommand {\smat}      [1] {\left[\begin{smallmatrix}{#1}}
\newcommand {\srix}          {\end{smallmatrix}\right]}
\newcommand {\s}      [1] {\begin{smallmatrix}{#1}}
\newcommand {\se}          {\end{smallmatrix}}

\newtheorem{defn}{Definition}[section]
\newtheorem*{remark}{Remark}

\newtheorem{kor}[defn]{Corollary}
\newtheorem{thm}[defn]{Theorem}

\title{Type II balanced truncation for deterministic bilinear control systems}

\author{Martin Redmann\thanks{Weierstrass Institute for Applied Analysis and Stochastics, Mohrenstr. 39, 10117 Berlin, Germany, Email: {\tt 
martin.redmann@wias-berlin.de}. Financial support by the DFG via Research Unit FOR 2402 is gratefully acknowledged.}}

\begin{document}

\maketitle

\begin{abstract}
When solving partial differential equations numerically, usually a high order spatial discretisation is needed.
Model order reduction (MOR) techniques are often used to reduce the order of spatially-discretised systems and hence reduce 
computational complexity. A particular MOR technique to obtain a reduced order model (ROM) is balanced truncation (BT), a method 
which has been extensively studied for deterministic linear systems. As so-called type I BT it has already been extended to bilinear equations, 
an important subclass of nonlinear systems. We provide an alternative generalisation of the linear setting to bilinear systems which is called type 
II BT. The Gramians that we propose in this context contain information about the control. It turns out that the new approach delivers energy bounds 
which are not just valid in a small neighbourhood of zero. Furthermore, we provide an $\mathcal{H}_\infty$-error bound which so far is not known when 
applying type I BT to bilinear systems.
\end{abstract}
\textbf{Keywords:} Model order reduction,  balanced truncation, bilinear systems, Gramians, error bound. 

\noindent\textbf{MSC classification:} 93A15, 93B05, 93B07, 93B36, 93C10.


\section{Introduction}

Numerical simulations are one of the conventional methods to study physical phenomena of dynamical systems. However, extracting all the complex 
system dynamics generally leads to large state-space systems, whose direct simulations are inefficient and involve huge computational 
cost. Hence, there is a need to consider model order reduction (MOR), aiming at replacing these large-scale systems by systems of much smaller state 
dimension. MOR for linear systems has been investigated intensively in recent years and is widely used in numerous applications, see e.g. 
\cite{antoulas, bennermehrmann, schildersvorstrommes}. In this work, we consider MOR for bilinear control systems, which can be considered as a 
bridge 
between linear and nonlinear systems. Applications of bilinear systems can be seen in various fields \cite{brunietal,Mohler,rugh}. \medskip

Several methods for linear systems have been extended to bilinear systems such as balanced truncation (BT) \cite{bennerdamm} or other balancing 
related methods \cite{hartmann}. Moreover, interpolation-based MOR has been applied \cite{bennerbreiten1, bennerbreiten2, breitendamm, flaggserkan}. 
In this manuscript, we focus on BT for bilinear systems which for linear systems has been studied in e.g. \cite{antoulas, moore}. Later on, the 
balancing concept for general nonlinear systems has been extended in a series of papers, see e.g. \cite{fujscherpen,grayscherpen,Scherpen}.\medskip

BT relies on controllability/reachability and observability Gramians. In \cite{bennerdamm} Gramians were proposed which we will call type I 
(bilinear) Gramians in this paper. The drawback of this approach is that only local energy bounds are available \cite{graymesko}. 
Furthermore, no error bound has been proved so far. The type I bilinear Gramians play a role for stochastic systems as well \cite{bennerdamm, 
redmannbenner,redSPA}, where they are also used in the context of balancing. Recently, a second 
way of balancing for stochastic systems was discussed \cite{bennerdammcruz,dammbennernewansatz,redmannspa2,BTtyp2EB}. It is based on another 
reachability Gramian and called type II ansatz. With this approach an $\mathcal H_\infty$-error bound can be achieved which cannot be proved in the 
ansatz used in \cite{redmannbenner, redSPA}. \medskip

In this paper, we introduce type II bilinear Gramians in Section \ref{sec:reach} which are inspired by the type II stochastic Gramians. The 
difference lies in additional information about the control in the bilinear Gramians. Under the assumption of having bounded controls, we then prove 
bounds for the reachability and observability energy of the underlying bilinear system using the type II bilinear Gramians. This provides a 
motivation to balance the bilinear system based on the new Gramians. The procedure will be explained in Section \ref{proctype2spa}. In Section 
\ref{spaerrorbound2}, an $\mathcal{H}_\infty$-error bound for type II bilinear BT will be proved, again assuming bounded controls. This error bound 
is the main result of this paper and has the same structure as in the linear case. Error bounds for BT applied to bilinear system did not exist 
before.

\section{Setting and type II Gramians}
\label{sec:reach}

We consider the following bilinear deterministic equation  
\begin{align}\label{stochstatenew}
             \frac{dx(t)}{dt}=Ax(t)+Bu(t)+\sum_{i=1}^m N_i x(t)u_i(t),\;\;\;x(0)=x_0, \;\;\;t\geq 0,
            \end{align}
where $A\in \mathbb R^{n\times n}$, $B\in \mathbb R^{n\times m}$ and $N_1, \ldots, N_m\in\mathbb R^{n\times n}$. Below, $x(t, x_0, u)$, $t\geq 
0$, denotes the solution to (\ref{stochstatenew}) with initial condition $x_0\in\mathbb R^n$ and control process $u=(u_1, \ldots, u_m)^T$. 
In our framework the state equation (\ref{stochstatenew}) is additionally equipped with an output equation \begin{align}\label{observeq2}
y(t, x_0, u)= C x(t, x_0, u),\;\;\;t\geq 0,
\end{align}
where $C\in\mathbb R^{p\times n}$. The control $u$ is usually assumed to be an $L^2_T$ function meaning that
     \begin{align*}
\left\|u\right\|_{L^2_T}^2:=\int_0^T u^T(t) u(t) dt=\int_0^T \left\|u(t)\right\|_2^2 dt<\infty
\end{align*}                            
for every $T>0$. Moreover, a classical assumption is the asymptotic stability of the uncontrolled equation (\ref{stochstatenew}), that is 
\begin{align}\label{assumstabasym}
\left\|x(t, x_0, 0)\right\|_2\rightarrow 0 
\end{align}
for $t\rightarrow \infty$ and every $x_0\in\mathbb R^n$, being equivalent to the Hurwitz property of $A$. Later on we will introduce a further 
condition on $u$ and a stronger assumption on the stability. This is required to define the type II Gramians. \smallskip

Before alternative Gramians are discussed, the existing theory will be summarised below. 

\paragraph{Type I reachability and observability Gramians}
 
In \cite{bennerdamm} and \cite{graymesko} bounds for the controllability/reachability and observability energy of the equations (\ref{stochstatenew}) 
and (\ref{observeq2}) have been proved using certain Gramians. We call these matrices type I Gramians here. The type I reachability Gramian $P_1$ is 
the unique solution to \begin{align}\label{gengenlyap}
A P_1+P_1 A^T+\sum_{i=1}^m N_i P_1 N_i^T = -B B^T,
\end{align}
whereas the type I observability Gramian is defined as the unique solution of
\begin{align}\label{gengenlyapobs}
A^T Q_1+Q_1 A+\sum_{i=1}^m N_i^T Q_1 N_i = -C^T C.
\end{align}
Equations (\ref{gengenlyap}) and (\ref{gengenlyapobs}) are also considered in the context of model order reduction for stochastic systems 
\cite{bennerdamm, redmannbenner}. A further discussion about these Gramians can be found in \cite{redmannspa2}. A unique positive semidefinite 
solution to the identities (\ref{gengenlyap}) and (\ref{gengenlyapobs}) exists if the following stability condition holds: 
\begin{align}\label{propmeansqstab}
 \sigma\left(I_n\otimes A+A\otimes I_n+\sum_{i=1}^m N_i\otimes N_i\right)\subset \mathbb C_-.
\end{align}
Property (\ref{propmeansqstab}) is also called asymptotic mean square stability because it represents a stability concept for stochastic 
systems \cite{redmannbenner,damm,staboriginal, redmannspa2}. It is satisfied if $A$ is asymptotically stable and the matrices $N_i$ 
are relatively 
small (in some norm) compared to the eigenvalues of $A$. That is why it is actually enough to assume (\ref{assumstabasym}) because equation 
(\ref{stochstatenew}) can be rescaled as \begin{align*}
 \frac{dx(t)}{dt}=Ax(t)+[\frac{1}{\gamma}B][\gamma u(t)]+\sum_{i=1}^m [\frac{1}{\gamma}N_i] x(t) [\gamma u_i(t)],
            \end{align*}
where the weighted matrices $\tilde N_i=\frac{1}{\gamma}N_i$ can be made arbitrary small with a sufficiently large constant $\gamma>0$.\medskip

No let us introduce the energy functionals. As in \cite{bennerdamm, graymesko, Scherpen} the controllability energy is
\begin{align*}
E_c(x_0):=\min_{{u\in L^2(]-\infty, 0]) \atop x(-\infty,x_0,u)=0}} \int_{-\infty}^0 \left\|u(t)\right\|_2^2 dt.
\end{align*}
In \cite{bennerdamm} and \cite{graymesko} the observability energy is considered for the case $B=0$, whereas for general non-linear systems $u\equiv 
0$ is assumed \cite{Scherpen}. With $B=0$, the observability energy is \begin{align}\label{eofuc}
E_o(x_0):=\max_{{u\in L^2([0, \infty[) \atop \left\|u\right\|_{L^2}<\alpha, B=0}} \int_{0}^\infty \left\|y(t, x_0, u)\right\|_2^2 dt,\;\;\;\text{for 
}\alpha>0\text{ small.}
\end{align}
The next theorem is a result from \cite{bennerdamm} and \cite{graymesko}.
\begin{thm}\label{bilinear energy}
Let $A$ be asymptotically stable (and at least one $N_i$ of full rank). If there are positive definite solutions $P_1$ and $Q_1$ to 
(\ref{gengenlyap}) and (\ref{gengenlyapobs}), then there exists a neighbourhood $\hat W(0)$ of zero such that \begin{align*}
      x_0^T P_1^{-1}x_0\leq E_c(x)\;\;\;\text{and}\;\;\;   E_o(x)\leq x_0^T Q_1x_0              \;\;\;\text{for }x\in \hat W(0).     
                                                                                      \end{align*}
\end{thm}
The above inequalities allow us to identify the hardly controllable and observable states (at least locally). Those are contained in the 
eigenspaces 
spanned by the eigenvectors of $P_1$ and $Q_1$, respectively, that correspond to the small eigenvalues. \\
More accurate bounds than in Theorem \ref{bilinear energy} were obtained for truncated Gramians \cite{bennergoyalredmann}, which are 
also computationally cheaper than the type I Gramians.

\paragraph{Type II reachability and observability Gramians}
We now introduce alternative Gramians which we will see, guarantee an error bound for the bilinear system. The so called type II Gramians 
are inspired by the stochastic case. A postive definite reachability Gramian $P_2$ which solves 
\begin{align}\label{newgram2}
 A^T P_2^{-1}+P_2^{-1}A+\sum_{i=1}^m N_i^T P_2^{-1} N_i \leq -P_2^{-1}BB^T P_2^{-1}
                                       \end{align} 
was initially introduced in \cite{dammbennernewansatz} in order to show the existence of an $\mathcal H_\infty$-error bound for BT applied to 
stochastic systems. There the balancing was based on $P_2$ and the type I Gramian $Q_1$, the solution to (\ref{gengenlyapobs}). The stochastic 
Gramian $P_2$ was furthermore analysed in \cite{bennerdammcruz} and used in \cite{redmannspa2,BTtyp2EB} in a more general form.  \medskip

An inequality is considered in (\ref{newgram2}) since the existence of a positive definite solution is not ensured when having an equality. A 
solution to inequality (\ref{newgram2}) exists if condition (\ref{propmeansqstab}) is satisfied \cite{dammbennernewansatz, redmannspa2}. As already 
mentioned 
above, (\ref{propmeansqstab}) can be weakened to the assumption of asymptotic stability of $A$ since $N_i$ can be made arbitrary small. \medskip

We will not take the matrices $Q_1$ and $P_2$ (stochastic type II balancing) to introduce a type II approach for bilinear systems. Moreover, we 
will modify them further. The idea is to let information about the control enter the new Gramians. This is done by choosing a constant $k>0$. So, 
the 
type II Gramians of the bilinear system (\ref{stochstatenew}), (\ref{observeq2}) are given by a perturbed matrix $A$: 
\begin{align}\label{bilinearreachgramtyp2}
 (A+\frac{k^2}{2} I_n)^T P^{-1}+P^{-1}(A+\frac{k^2}{2}I_n)+\sum_{i=1}^m N_i^T P^{-1} N_i &\leq -P^{-1}BB^T P^{-1},\\ \label{bilinearobservgramtyp2}
 (A+\frac{k^2}{2} I_n)^T Q+Q (A+\frac{k^2}{2} I_n)+\sum_{i=1}^m N_i^T Q N_i &= -C^T C.
       \end{align}
At the same time we suppose to have a control which is uniformly bounded by the perturbing constant on a finite time interval $[0, T]$ 
\begin{align}\label{boundforu}
               \left\|u(t)\right\|_2\leq k,\;\;\;t\in[0, T].
                                                                                                                            \end{align}
Of course $k$ cannot be arbitrary large because we need the existence of the solutions to (\ref{bilinearreachgramtyp2}) and 
(\ref{bilinearobservgramtyp2}). Ideally $k$ is such that \begin{align}\label{propmeansqstabpert}
 \sigma\left(I_n\otimes (A+\frac{k^2}{2} I_n)+(A+\frac{k^2}{2} I_n)\otimes I_n+\sum_{i=1}^m N_i\otimes N_i\right)\subset \mathbb C_-
\end{align}
holds. We observe that the type II Gramians $P, Q$ are coupled with the control $u$. We allow large controls ($L^\infty$-sense) if the system is 
relatively stable (largest real part of the eigenvalues of $A$ is small and $N_i$ are small) and only small controls are admissible if the system is 
close to be unstable. We can again weaken condition (\ref{propmeansqstabpert}) and only assume  $\sigma\left(A+\frac{k^2}{2}\right)\subset \mathbb 
C_-$ 
since we can rescale $N_i$, i.e., we replace it by $\tilde N_i=\frac{1}{\gamma}N_i$. In this situation we pay a price for the rescaling since if we 
want to bound the rescaled control $\tilde u=\gamma u$ as in (\ref{boundforu}), it is required to have $ \left\|u(t)\right\|_2\leq \frac{k}{\gamma}$ 
which restricts the controls even more. From now on, we assume that (\ref{boundforu}) and (\ref{propmeansqstabpert}) hold at the same time, knowing 
that we need a less restrictive assumption than (\ref{propmeansqstabpert}) when using a smaller bound for the controls.\medskip

Let us now investigate how much energy is necessary if we control the system from zero into a certain direction. We desire to bound the reachability 
energy with the type II Gramian $P$. Let $(p_{k})_{k=1, \ldots, n}$ be eigenvectors of $P$ such that they represent an orthonormal basis of $\mathbb 
R^n$. The corresponding eigenvalues are denoted by $(\lambda_{k})_{k=1, \ldots, n}$. Then, for $t\in[0, T]$
\begin{align*}
\langle x(t, 0, u), p_{k}  \rangle_2^2 &\leq \lambda_{k}\; \sum_{i=1}^n \lambda_{i}^{-1} \langle x(t, 0, u), p_{i}  
\rangle_2^2=\lambda_{k} \;\left\|\sum_{i=1}^n \lambda_{i}^{-\frac{1}{2}} \langle x(t, 0, u), p_{i}  \rangle_2 \;p_{i}\right\|_{2}^2
\\& =\lambda_{k}\; \left\|P^{-\frac{1}{2}} x(t, 0, u)\right\|_{2}^2 = \lambda_{k} \;\left[x(t, 0, u)^T P^{-1} x(t, 0, u)\right].
\end{align*}
To shorten the notation we write $x(t)$ instead of $x(t, 0, u)$ from time to time below. So, by the product rule and by inserting 
(\ref{stochstatenew}), we have \begin{align}\label{prodruleP}
x^T(t)P^{-1} x(t)=&\int_0^t x^T(s)P^{-1} dx(s)+\int_0^t dx^T(s) P^{-1} x(s)\\ \nonumber=&\int_0^t x^T(s) P^{-1} A x(s) ds+\int_0^t x^T(s)P^{-1} B 
u(s) ds\\ \nonumber
&+\int_0^t x^T(s) A^T P^{-1} x(s)ds+\int_0^t u^T(s)B^T P^{-1} x(s)ds\\ \nonumber
&+\sum_{i=1}^m \left(\int_0^t x^T(s) P^{-1} N_i x(s) u_i(s) ds+\int_0^t x^T(s) u_i(s)N_i^T P^{-1} x(s)ds\right).
\end{align}
We analyse the last term above which can be written as \begin{align*}
&\sum_{i=1}^m 2 \int_0^t x^T(s) P^{-1} N_i x(s) u_i(s) ds=\sum_{i=1}^m 2 \int_0^t \left\langle  P^{-\frac{1}{2}} x(s)u_i(s), P^{-\frac{1}{2}} N_i 
x(s) \right\rangle_2 ds\\&\leq \sum_{i=1}^m \left( \int_0^t \left\|  P^{-\frac{1}{2}} x(s)u_i(s)\right\|_2^2ds + \int_0^t\left\|P^{-\frac{1}{2}} N_i 
x(s) \right\|^2_2 ds\right)\\&=\int_0^t  x^T(s) P^{-1} x(s) \left\|u(s)\right\|_2^2 ds + \sum_{i=1}^m \int_0^t x^T(s)N_i^TP^{-1} N_i x(s)ds\\&\leq 
\int_0^t  
x^T(s) P^{-1}k^2 x(s) ds + \sum_{i=1}^m \int_0^t x^T(s)N_i^TP^{-1} N_i x(s)ds
\end{align*}
using the bound in (\ref{boundforu}). Summarising the above steps, we obtain 
\begin{align*}
\langle x(t, 0, u), p_{k}  \rangle_2^2 &\leq \lambda_{k}\left[\int_0^t x^T(s)(A^T P^{-1}+P^{-1}A+ \sum_{i=1}^m N_i^T P^{-1} N_i+k^2 P^{-1}) 
x(s)\right. ds\\ &\ \ \ \ \ \ \ \ \ \ \ \ \ \ \left.+2\int_0^t x^T(s) P^{-1} B u(s) ds\right].
\end{align*}
Plugging in (\ref{bilinearreachgramtyp2}) yields  \begin{align*}
 \langle x(t, 0, u), p_{k}  \rangle_2^2 &\leq \lambda_{k}\int_0^t 2 x^T(s) P^{-1} B u(s)-x^T(s) P^{-1}BB^T P^{-1} x(s)ds\\
&=\lambda_{k} \int_0^t \left\|u(s)\right\|_2^2 ds - \left\|B^T P^{-1} x(s)-u(s)\right\|_2^2ds.
\end{align*}
Consequently, we have \begin{align}\label{diffreachjaneintype2}
\lambda_{k}^{-\frac{1}{2}} \sup_{t\in[0, T]}\left\vert\langle x(t, 0, u), p_{k}  \rangle_2\right\vert \leq \left\|u\right\|_{L^2_T}.
\end{align}
So, by (\ref{diffreachjaneintype2}), large controllability energy is needed if $\lambda_{k}$ is small especially when a large component in the 
direction of $p_{k}$ shall be reached (a component on or outside the unit sphere). This implies that difficult to reach states have a ``large'' 
component in the eigenspaces of $P$ belonging to the small eigenvalues.
\begin{remark}
We can replace $P$ by the stochastic type II Gramian $P_2$ which satisfies (\ref{newgram2}). This results in the following inequality \begin{align*}
x^T(t)P_2^{-1}x(t) &\leq \left[\int_0^t x^T(s)(A^T P_2^{-1}+P_2^{-1}A+ \sum_{i=1}^m N_i^T P_2^{-1} N_i) 
x(s)\right. ds\\ &\ \ \ \ \ \ \ \ \ \left.+2\int_0^t x^T(s) P_2^{-1} B u(s) ds+\int_0^t x^T(s)P_2^{-1}x(s) \left\|u(s)\right\|_2^2 ds\right].
\end{align*}
Inserting (\ref{newgram2}), we then see that \begin{align*}
 x^T(t)P_2^{-1}x(t) \leq \int_0^t \left\|u(s)\right\|_2^2ds +\int_0^t x^T(s)P_2^{-1}x(s) \left\|u(s)\right\|_2^2 ds .
\end{align*}
By Gronwall's inequality, we obtain \begin{align*}
 x^T(t)P_2^{-1}x(t) \leq \int_0^t \left\|u(s)\right\|_2^2ds \expn^{\int_0^t \left\|u(s)\right\|_2^2 ds}
\end{align*}
which leads to inequality (\ref{diffreachjaneintype2}) with an additional exponential term but in this case no bound on the control is assumed. 
Let us suppose that $\left\|u\right\|_{L^2_T}\leq 1$ in case $P_2$ is used. Then, a small eigenvalue $\lambda_{2, k}$ of $P_2$ implies that $\langle 
x(t, 0, u), p_{2, k}  \rangle_2$ ($p_{2, k}$ is the corresponding eigenvector) is small. This means that no large component in the direction of 
$p_{2, k}$ can be reached with a small control.
\end{remark}
Let us now turn our attention to the type II Gramian $Q$. We shorten the notation again and write $x_{x_0}(t)$ instead of $x(t, x_0, u)$. The product 
rule 
yields \begin{align*}
x_{x_0}^T(t)Q x_{x_0}(t)-x_0^TQx_0=&\int_0^t x_{x_0}^T(s)Q dx_{x_0}(s)+\int_0^t dx_{x_0}^T(s) Q x_{x_0}(s),
\end{align*}
where $t\in[0, T]$. Following the steps from (\ref{prodruleP}) onwards, we obtain \begin{align*}
x_{x_0}^T(t)Q x_{x_0}(t)-x_0^TQx_0 &\leq \left[\int_0^t x_{x_0}^T(s)(A^T Q+Q A+ \sum_{i=1}^m N_i^T Q N_i+k^2 Q) 
x_{x_0}(s)\right. ds\\ &\ \ \ \ \ \ \ \ \left.+2\int_0^t x_{x_0}^T(s) Q B u(s) ds\right].
\end{align*}
We insert (\ref{bilinearobservgramtyp2}) and evaluate the functions at the final time which gives \begin{align}\label{observattyp2en}
\int_0^T \left\|y(s, x_0, u)\right\|_2^2 ds\leq x_0^TQx_0+2\int_0^T x(s, x_0, u)^T Q B u(s) ds.
\end{align}
As in (\ref{eofuc}) we assume $B=0$. This is a natural choice since in the observability problem an unknown initial condition shall be reconstructed 
from the observations $y$. Since $Bu$ is a term that does not depend on $x_0$, it can be assumed to be known and hence be neglected in the 
considerations. With $B=0$ and (\ref{observattyp2en}) we see that the states which produce little observation 
energy are close to the kernel of $Q$. They are contained in the eigenspaces of $Q$ corresponding to the small eigenvalues.
\begin{remark}
If we use $Q_1$ which satisfies (\ref{newgram2}) instead of $Q$, we get an extra term in the energy bound but in this case there is no bound on $u$. 
So, we have \begin{align*}
\int_0^T \left\|y(s, x_0, u)\right\|_2^2 ds\leq &x_0^TQ_1x_0+2\int_0^T x(s, x_0, u)^T Q_1 B u(s) ds\\
&+\int_0^T x(s, x_0, u)^T Q_1 x(s, x_0, u) \left\|u(s)\right\|_2^2 ds.
\end{align*}
Now, if we say that the control $u$ is small, we can conclude that the observation energy is small if the initial state is close to the kernel of 
$Q_1$.
\end{remark}

\section{Type II balanced truncation}
\label{proctype2spa}

Before considering an $\mathcal{H}_\infty$- error bound for BT based on the type II Gramians $P$ and $Q$, we summarise the theory of balancing which 
is similar to the deterministic linear case \cite{antoulas, moore}. 

States that are difficult to reach can be characterised by $P$, cf. (\ref{diffreachjaneintype2}). These states have large components in the span of 
the eigenvectors corresponding to small eigenvalues of the reachability Gramian $P$. Similarly, states that are difficult to 
observe are the ones that have large components in the span of eigenvectors corresponding to small eigenvalues of the observability Gramian $Q$, see 
(\ref{observattyp2en}). Now, balancing a system relies on the idea to create a system, where dominant reachable and 
observable states are the same, i.e., reachability and observability Gramians are simultaneously transformed such that they are equal and diagonal. 
BT for bilinear systems based on the Gramians $P_1$ and $Q_1$ (type I ansatz) has already been studied \cite{bennerdamm} and for related energy 
functionals compare \cite{graymesko}. For type I BT no error bound could be shown so far. Now, the procedure for the type II approach is explained. 
This ansatz allows us to show an $\mathcal H_\infty$-error bound in Section \ref{spaerrorbound2}.\medskip

We consider a control system consisting of state equation (\ref{stochstatenew}) and output equation (\ref{observeq2}) 
\begin{align}\label{controlsys2}
\frac{d{x}(t)}{dt}&={A} x(t) + B u(t)+\sum_{i=1}^m N_i x(t) u_i(t),\\ \nonumber
y(t)&= {C} x(t),\;\;\;t\geq 0,
\end{align}
Recall that the state equation in (\ref{controlsys2}) is asymptotically stable, i.e., $\sigma\left(A+\frac{k^2}{2}\right)\subset \mathbb C_-$ or 
ideally (\ref{propmeansqstabpert}) is satisfied.
Introducing a transformation matrix $T\in\mathbb{R}^{n\times n}$ which is assumed to be non-singular, the states are transformed as follows: 
\[\hat{x}(t) = Tx(t),\] 
such that system (\ref{controlsys2}) becomes
\begin{align}\label{balancingtransformation}
\frac{d}{dt}\hat{x}(t)&=\hat{A} \hat{x}(t) + \hat{B}u(t)+\sum_{i=1}^m \hat N_i x(t) u_i(t),\\ \nonumber
y(t)&= \hat{C} \hat{x}(t),\;\;\;t\geq 0,
\end{align}
where $\hat{A} = TAT^{-1}$, $\hat{B}=TB$, $\hat{C} =CT^{-1}$ and $\hat{N_i} = TN_i T^{-1}$. The input-output 
map remains the same, only the state and the systems matrices are transformed.  

$P$ and $Q$, the reachability and observability Gramians of system (\ref{controlsys2}),
which satisfy (\ref{bilinearreachgramtyp2}) and (\ref{bilinearobservgramtyp2}) can be transformed into reachability and observability Gramians of the 
transformed system (\ref{balancingtransformation}): \begin{align*}
\hat{P} =TPT^T\;\;\; \text{and}\;\;\; \hat{Q} = T^{-T}QT^{-1}.
\end{align*}
The above relation is obtained by multiplying (\ref{bilinearreachgramtyp2}) and (\ref{bilinearobservgramtyp2}) with 
$T^{-T}$ from the left and $T^{-1}$ from the right. The Hankel singular values (HSVs) $\sigma_1\ge\ldots\ge\sigma_n$, where $\sigma_i 
=\sqrt{\lambda_i(P Q)}$ ($i=1,\ldots,n$), of the original and transformed system are the same. The above transformation is a balancing 
transformation if the transformed Gramians are equal and diagonal. Such a transformation 
always exists if $Q>0$ (observation energy is always non zero for every $x_0\neq 0$). We also need that $P>0$ but this is automatically satisfied. A 
balanced system is obtained by choosing 
\[
T=\Sigma^{-\frac{1}{2}}U^T L^T \quad\text{and}\quad T^{-1}=KV\Sigma^{-\frac{1}{2}},
\]
where $\Sigma=\diag(\sigma_{1},\ldots,\sigma_n)>0$ is the diagonal matrix of HSVs. $Y$, $Z$, $L$ and $K$ are computed as follows. Let $P = KK^T$, 
$Q=LL^T$ be square root factorisations of $P$ and $Q$, then an SVD of $K^TL = V\Sigma U^T$ gives the required matrices. With this transformation 
$\hat{P}=\hat{Q}=\Sigma$. This implies that $\Sigma$ characterises both the reachability and observability in system 
(\ref{balancingtransformation}). The smaller the diagonal entry of $\Sigma$, the less important the corresponding state component in the system 
dynamics of (\ref{balancingtransformation}).

Below, let $T$ be the balancing transformation as stated above, then we partition the coefficients of the balanced realisation as 
follows:\begin{align*}
T{A}T^{-1}= \smat{A}_{11}&{A}_{12}\\ 
{A}_{21}&{A}_{22}\srix,\;\;\; T{B} = \smat{B}_1\\ {B}_2\srix,\;\;\;  
{CT^{-1}} = \smat{C}_1 &{C}_2\srix,\;\;\;T{N_i}T^{-1}= \smat{N}_{i, 11}&{N}_{i, 12}\\ 
{N}_{i, 21}&{N}_{i, 22}\srix, \end{align*}
where ${A}_{11}\in\R^{r\times r}$ etc. Furthermore, by setting $\hat x= \smat{x}_1\\ 
{x}_2\srix$, where $x_1(t)\in\mathbb{R}^r$, we obtain the transformed partitioned system 
{\footnotesize\begin{align} 
\frac{d}{dt}\mat{c} {x}_1(t)\\{x}_2(t)\rix \hspace{-0.05cm}&=\hspace{-0.05cm}\mat{cc} {A}_{11}&\hspace{-0.25cm}{A}_{12}\\ 
{A}_{21}&\hspace{-0.25cm}{A}_{22}\rix\hspace{-0.1cm} \mat{c} {x}_1(t)\\ {x}_2(t)\rix + \mat{c} 
{B}_1\\ {B}_2\rix \hspace{-0.05cm} u(t)\hspace{-0.05cm}+\sum_{i=1}^m \hspace{-0.05cm}\mat{cc} {N}_{i, 11}&\hspace{-0.25cm}{N}_{i, 12}\\ 
{N}_{i, 21}&\hspace{-0.25cm}{N}_{i, 22}\rix \hspace{-0.1cm}\mat{c} {x}_1(t)\\ {x}_2(t)\rix \hspace{-0.05cm} u_i(t),\label{balrelpart}\\ 
y(t)&= \mat{cc} {C}_1 &{C}_2\rix \mat{c} {x}_1(t)\\ {x}_2(t)\rix,\;\;\;t\geq 0.\label{balrelpartout}
\end{align}}
\hspace{-0.05cm}From this system we aim to obtain a approximating system with reduced dimension $r\ll n$. For BT the second row in (\ref{balrelpart}) 
is truncated and the remaining $x_2$ components in the first row of (\ref{balrelpart}) and in (\ref{balrelpartout}) are set to zero. This leads to a 
ROM having the same structure as (\ref{controlsys2}), that is,
\begin{align}\label{generalreducedsys}
\frac{dx_r(t)}{dt}&=A_{11} x_r(t)+ B_1 u(t)+\sum_{i=1}^m {N}_{i, 11} x_r(t) u_i(t),\\ \nonumber
y_r(t)&= C_1 x_r(t),\;\;\;t\geq 0,
\end{align}
where $A_{11}, N_{i, 11}\in\mathbb R^{r\times r}$, $B_1 \in\mathbb R^{r\times m}$ and $C_1\in\mathbb R^{p\times r}$.
In equations (\ref{balrelpart}) and (\ref{balrelpartout}), the difficult to reach and observe states are 
represented by $x_2$, which correspond to the smallest HSVs  $\sigma_{r+1}, \ldots, \sigma_n$, but of course $r$ has to be chosen such that the 
neglected HSVs are small ($\sigma_{r+1}\ll\sigma_{r}$). 

\section{$\mathcal H_\infty$-error bound for type II BT}\label{spaerrorbound2}

In this section, we show that the $\mathcal H_\infty$-error bound, that is known from the linear case \cite{antoulas}, is also true for bilinear 
systems when using the type II approach. Unfortunately, this result could no yet be established when using the type I Gramians.\medskip

We recall the original model that we aim to reduce:
\begin{align}\label{stochstatenewextend}
             \frac{dx(t)}{dt}&=Ax(t)+Bu(t)+\sum_{i=1}^m N_i x(t)u_i(t),\;\;\;x(0)=0,\\ \nonumber
             y(t)&=Cx(t), \;\;\;t\geq 0,
            \end{align}
where the matrices and vectors above are partitioned as follows \begin{align*}
 A=\smat{A}_{11}&{A}_{12}\\ 
{A}_{21}&{A}_{22}\srix,\;x=\smat x_1 \\ x_2\srix,\;B=\smat B_1 \\ B_2\srix,\; N_i=\smat{N}_{i, 11}&{N}_{i, 12}\\ 
{N}_{i, 21}&{N}_{i, 22}\srix,\;C= \smat C_1 & C_2\srix.
            \end{align*}
Below, we prove an $\mathcal H_\infty$-error bound when the balancing relies on the matrix (in)equalities (\ref{bilinearreachgramtyp2}) and 
(\ref{bilinearobservgramtyp2}). We replace equation (\ref{bilinearobservgramtyp2}) by an inequality because we do not need the equality in the 
proof. To simplify the notation, we assume that system (\ref{stochstatenewextend}) is balanced already, i.e., we have applied the balancing 
transformation from Section \ref{proctype2spa} already. Hence, the Gramians $P$ and $Q$ are equal and coincide with the diagonal matrix 
$\Sigma=\diag(\Sigma_1, \Sigma_2)$, where $\Sigma_1=\diag(\sigma_1, \ldots, \sigma_r)$ is the matrix of large and $\Sigma_2=\diag(\sigma_{r+1}, 
\ldots, \sigma_n)$ the matrix of neglected small HSVs. \\
The following ROM is supposed to be compared with the original model (\ref{stochstatenewextend}):
\begin{align}\label{romstochstatenewextend}
             \frac{dx_r(t)}{dt}&=A_{11}x_r(t)+B_1u(t)+\sum_{i=1}^m N_{i, 11} x_r(t)u_i(t),\;\;\;x_r(0)=0,\\ 
\nonumber     y_r(t)&=C_1x_r(t), \;\;\;t\geq 0.
            \end{align}
The next theorem deals with the $L_T^2$-error between the full and the reduced order output.       
\begin{thm}\label{thmfirsthsv}
Let $x(0)=0$, $x_{r}(0)=0$ and $\left\|u(t)\right\|_2\leq k$, $t\in[0, T]$, $T>0$, where $k$ is the constant that enters in 
(\ref{bilinearreachgramtyp2}) and (\ref{bilinearobservgramtyp2}). Then,
\begin{align*}
   \left\|y-y_r\right\|_{L^2_T}\leq 2 (\tilde \sigma_{1}+\tilde \sigma_{2}+\ldots + \tilde \sigma_\nu)  \left\|u\right\|_{L^2_T},            
  \end{align*}
where $y$ is the output of the original system (\ref{stochstatenewextend}), $y_r$ is the output of the type II BT approach ROM and 
$\tilde\sigma_{1}, \tilde\sigma_{2}, \ldots,\tilde\sigma_\nu$ are the distinct diagonal entries of 
$\Sigma_2=\diag(\sigma_{r+1},\ldots,\sigma_n)=\diag(\tilde\sigma_{1} I, \tilde\sigma_{2} I, \ldots, \tilde\sigma_\nu I)$.
\begin{proof}
We sometimes omit the time dependence of the functions in this proof to keep the notation as easy as possible. Inserting for $y$ and $y_r$ yields
 \begin{align*} -\left\|y-y_r\right\|^2_2=-\left\|C_1(x_1-x_r)+C_2 x_2\right\|^2_2=-\mat{c}x_1-x_r \\ x_2\rix^T C^T C \mat{c}x_1-x_r \\ x_2\rix.
\end{align*}
The partitioned matrix (in)equality (\ref{bilinearobservgramtyp2}) \begin{align}\nonumber
  &\begin{smallmatrix}\smat{A}_{11}&{A}_{12}\\ 
{A}_{21}&{A}_{22}\srix^T \smat{\Sigma}_{1}&\\ 
&{\Sigma}_{2}\srix+\smat{\Sigma}_{1}&\\ 
&{\Sigma}_{2}\srix \smat{A}_{11}&{A}_{12}\\ 
{A}_{21}&{A}_{22}\srix +\end{smallmatrix}\sum_{i=1}^m\begin{smallmatrix}\smat{N}_{i, 11}&{N}_{i, 12}\\ 
{N}_{i, 21}&{N}_{i, 22}\srix^T \smat{\Sigma}_{1}&\\ 
&{\Sigma}_{2}\srix \smat{N}_{i, 11}&{N}_{i, 12}\\ 
{N}_{i, 21}&{N}_{i, 22}\srix+k^2\smat{\Sigma}_{1}&\\ 
&{\Sigma}_{2}\srix\end{smallmatrix}\\ \label{observeineq} &\begin{smallmatrix}\leq {-C^T C}\end{smallmatrix}
                                       \end{align}
leads to  \begin{align*} &-\left\|y-y_r\right\|^2_2\geq\\ & 2 (x_1-x_r)^T \Sigma_1  \smat {A}_{11}&{A}_{12} 
\srix \smat x_1-x_r \\ x_2\srix+\sum_{i=1}^m\left(\smat {N}_{i, 11}&{N}_{i, 12} \srix \smat x_1-x_r \\ 
x_2\srix\right)^T \Sigma_1 \smat {N}_{i, 11}&{N}_{i, 12} \srix \smat x_1-x_r \\ 
x_2\srix\\ &+(x_1-x_r)^T \Sigma_1 k^2 (x_1-x_r)+x_2^T \Sigma_2 k^2  x_2 \\ & +2 x_2^T 
\Sigma_2  \smat {A}_{21}&{A}_{22} \srix 
\smat x_1-x_r \\ x_2\srix +\sum_{i=1}^m\left(\smat {N}_{i, 21}&{N}_{i, 22} \srix \smat x_1-x_r \\ 
x_2\srix\right)^T \Sigma_2 \smat {N}_{i, 21}&{N}_{i, 22} \srix \smat x_1-x_r \\ x_2\srix.
\end{align*}  
Using the above summands, we define {\allowdisplaybreaks
\begin{align*} \mathcal T_1:&= 2 (x_1-x_r)^T \Sigma_1 \smat {A}_{11}&{A}_{12} 
\srix \smat x_1-x_r \\ x_2\srix,\\ \mathcal T_2:&= \sum_{i=1}^m\left(\smat {N}_{i, 11}&{N}_{i, 12} \srix \smat x_1-x_r \\ 
x_2\srix\right)^T \Sigma_1 \smat {N}_{i, 11}&{N}_{i, 12} \srix \smat x_1-x_r \\ 
x_2\srix +(x_1-x_r)^T \Sigma_1 k^2 (x_1-x_r),\\
\mathcal T_3:&=2 x_2^T \Sigma_2  \smat {A}_{21}&{A}_{22} \srix \smat x_1-x_r \\ 
x_2\srix,\\ \mathcal T_4:&= \sum_{i=1}^m \left(\smat {N}_{i, 21}&{N}_{i, 22} \srix \smat x_1-x_r \\ 
x_2\srix\right)^T \Sigma_2 \smat {N}_{i, 21}&{N}_{i, 22} \srix \smat x_1-x_r \\ x_2\srix+ x_2^T \Sigma_2 k^2 x_2 .
\end{align*}}
The differential of $x_1-x_r$ is given by\begin{align*}
     \frac{d(x_1-x_r)}{dt}=\mat{cc} A_{11}&A_{12}\rix \mat{c} x_1-x_r \\ x_2 \rix +\sum_{i=1}^m \mat{cc} N_{i, 11}&N_{i, 12} \rix \mat{c} 
x_1-x_r\\ x_2 \rix u_i           
                \end{align*}
which we insert into the following definition: \begin{align*}
     D_1:=\frac{d\left((x_1(t)-x_r(t))^T\Sigma_1 (x_1(t)-x_r(t))\right)}{dt}=2(x_1(t)-x_r(t))^T\Sigma_1 \frac{d(x_1(t)-x_r(t))}{dt}.
                \end{align*}
Hence, we have \begin{align*}
     D_1=2(x_1-x_r)^T\Sigma_1 \left(\smat A_{11}&A_{12}\srix \smat x_1-x_r \\ x_2 \srix +\sum_{i=1}^m \smat N_{i, 11}&N_{i, 12} \srix 
\smat x_1-x_r\\ x_2 \srix u_i \right).
                \end{align*}
We use an elementary estimate for the following inner product \begin{align*}
     &\sum_{i=1}^m 2 (x_1-x_r)^T\Sigma_1 \smat N_{i, 11}&N_{i, 12} \srix \smat x_1-x_r\\ x_2 \srix u_i\\&=\sum_{i=1}^m 2 \left\langle 
\Sigma_1^{\frac{1}{2}}(x_1-x_r) u_i, \Sigma_1^{\frac{1}{2}} \smat N_{i, 11}&N_{i, 12} \srix \smat x_1-x_r\\ x_2 \srix \right\rangle_2\\&\leq 
\sum_{i=1}^m \left( \left\|\Sigma_1^{\frac{1}{2}}(x_1-x_r) u_i\right\|^2_2+\left\| \Sigma_1^{\frac{1}{2}} \smat N_{i, 11}&N_{i, 12} \srix \smat
x_1-x_r\\ x_2 \srix \right\|^2_2\right)\\&= (x_1-x_r)^T\Sigma_1(x_1-x_r) \left\|u\right\|^2_2+ \sum_{i=1}^m \left(\smat N_{i, 11}&N_{i, 
12} \srix \smat x_1-x_r\\ x_2 \srix \right)^T \Sigma_1\smat N_{i, 11}&N_{i, 12} \srix \smat x_1-x_r\\ x_2 \srix.
                \end{align*}
Using the fact that $u$ is bounded, i.e., $\left\|u\right\|_2\leq k$ then yields \begin{align*}
\frac{d}{dt} (x_1(t)-x_r(t))^T\Sigma_1 (x_1(t)-x_r(t))\leq\mathcal T_1+\mathcal T_2. 
\end{align*}         
From (\ref{stochstatenewextend}) the variable $x_2$ satifies \begin{align*}
     \frac{dx_2(t)}{dt}= \mat{cc}A_{21}& A_{22}\rix \mat{c}x_1(t) \\ x_2(t)\rix + B_2u(t) + \sum_{i=1}^m \mat{cc}N_{i, 21}& N_{i,22}\rix \mat{c} 
x_1(t) \\ x_2(t)\rix u_i(t).          
                \end{align*}
Plugging this into $\frac{d x_2(t)^T\Sigma_2 x_2(t)}{dt}=2x_2(t)^T\Sigma_2 \frac{x_2(t)}{dt}$ provides the following relation 
\begin{align}\label{derix2}
\frac{dx_2}{dt}=2x_2^T\Sigma_2 \left(\mat{cc}A_{21}& A_{22}\rix 
\mat{c}x_1 \\ x_2\rix + B_2u + \sum_{i=1}^m \mat{cc}N_{i, 21}& N_{i,22}\rix \mat{c} x_1 \\ x_2\rix u_i\right).
                   \end{align}
We find an upper bound for the last term the same way we have done before\begin{align*}
&\sum_{i=1}^m 2x_2^T\Sigma_2 \smat N_{i, 21}& N_{i,22}\srix \smat x_1 \\ x_2\srix u_i=\sum_{i=1}^m 2 \left\langle 
\Sigma_2^{\frac{1}{2}}x_2 u_i, \Sigma_2^{\frac{1}{2}} \smat N_{i, 21}&N_{i, 22} \srix \smat x_1\\ x_2 \srix \right\rangle_2\\ &\leq  
x_2^T\Sigma_2 x_2 k^2+ \sum_{i=1}^m \left(\smat N_{i, 21}&N_{i, 22} \srix \smat x_1\\ x_2 \srix \right)^T \Sigma_2\smat N_{i, 21}&N_{i, 22} 
\srix \smat x_1\\ x_2 \srix
                \end{align*}
applying again that the control is bounded by $k$. This yields \begin{align*}
     \frac{d}{dt}\left(x_2(t)^T\Sigma_2 x_2(t)\right)\leq&\left[\mathcal T_3+2x_2^T\Sigma_2 ({A}_{21}x_r+B_2 u)\right]\\
&+\left[\mathcal T_4+2 \sum_{i=1}^m (N_{i, 21} x_r)^T \Sigma_2 \mat{cc} N_{i, 21} & N_{i, 22}\rix \mat{c} x_1 \\x_2\rix \right.\\&\ \ 
\ \ \ \ \left.- \sum_{i=1}^m(N_{i, 21} x_r)^T \Sigma_2 N_{i, 21} x_r\right].
                \end{align*} 
Summarising the above computations, we obtain  \begin{align} \label{key1}
-\left\|y-y_r\right\|^2_2&\geq \hspace{-0.1cm} \mathcal T_1+\mathcal T_2+\mathcal T_3+\mathcal T_4\geq \hspace{-0.1cm} 
\frac{d}{dt}((x_1-x_r)^T\Sigma_1 (x_1-x_r))\hspace{-0.07cm}+\hspace{-0.07cm}\frac{d}{dt}(x_2^T\Sigma_2 x_2)\\&\nonumber -2x_2^T\Sigma_2 
({A}_{21}x_r+B_2 u)-2 \sum_{i=1}^m (N_{i, 21} x_r)^T \Sigma_2 \mat{cc} N_{i, 21} & N_{i, 22}\rix \mat{c} x_1 \\x_2\rix .                              
                                                              \end{align}
For the moment we assume that $\Sigma_2=\sigma I$. Since we have zero initial conditions, it holds that \begin{align}\label{firstestimate}
\int_0^T \left\|y(t)-y_r(t)\right\|^2_2 dt\leq  2\sigma^2 &\left(\int_0^T x_2^T\Sigma_2^{-1} ({A}_{21}x_r+B_2 u)dt\right.\\&\ \ \ +\left.\int_0^T 
\sum_{i=1}^m (N_{i, 21} x_r)^T \Sigma_2^{-1} \mat{cc} N_{i, 21} & N_{i, 22}\rix \mat{c} x_1 \\x_2\rix dt\right).  \nonumber                      
                                               \end{align}
Inequality (\ref{bilinearreachgramtyp2}) and the Schur complement condition on definiteness imply\begin{align}\label{schurposdef}
 \mat{cc} A^T \Sigma^{-1}+\Sigma^{-1}A+\sum_{i=1}^m N_i^T \Sigma^{-1} N_i +\Sigma^{-1}k^2 & \Sigma^{-1}B\\
 B^T \Sigma^{-1}& -I\rix\leq 0.
                                       \end{align}
If we multiply the matrix inequality (\ref{schurposdef}) with $\smat x_1+x_r \\ x_2 \\ 2u\srix^T$ from the left and with 
$\smat x_1+x_r \\ x_2\\ 2u\srix$ from the right, we get \begin{align*}
 4\left\|u\right\|^2_2\geq &\; 2 (x_1+x_r)^T \Sigma_1^{-1} \left( \smat {A}_{11}&{A}_{12} \srix \smat x_1+x_r \\ x_2\srix +2 B_1 
u\right)+(x_1-x_r)^T \Sigma_1^{-1}k^2(x_1-x_r) \\&+\sum_{i=1}^m\left(\smat {N}_{i, 11}&{N}_{i, 12} \srix \smat x_1+x_r \\ x_2\srix\right)^T 
\Sigma_1^{-1} \smat {N}_{i, 11}&{N}_{i, 12} \srix \smat x_1+x_r \\ x_2\srix \\ & +2  x_2^T \Sigma_2^{-1} \left( \smat {A}_{21}&{A}_{22} \srix 
\smat x_1+x_r \\ x_2\srix +2 B_2 u\right)+x_2^T \Sigma_2^{-1}k^2x_2\\& +\sum_{i=1}^m\left(\smat {N}_{i, 21}&{N}_{i, 22} \srix \smat x_1+x_r \\ 
x_2\srix\right)^T \Sigma_2^{-1} \smat {N}_{i, 21}&{N}_{i, 22} \srix \smat x_1+x_r \\ x_2\srix.
            \end{align*}
The above terms are used to define {\allowdisplaybreaks
\begin{align*} \mathcal T_5:&= 2 (x_1+x_r)^T \Sigma_1^{-1} \left( \smat {A}_{11}&{A}_{12} \srix \smat x_1+x_r \\ x_2\srix +2 B_1 
u\right),\\
\mathcal T_6:&=(x_1+x_r)^T \Sigma_1^{-1}k^2(x_1+x_r)+\sum_{i=1}^m\left(\smat {N}_{i, 11}&{N}_{i, 12} \srix \smat x_1+x_r \\ x_2\srix\right)^T 
\Sigma_1^{-1} \smat {N}_{i, 11}&{N}_{i, 12} \srix \smat x_1+x_r \\ x_2\srix, \\
\mathcal T_7:&=2  x_2^T \Sigma_2^{-1} \left( \smat {A}_{21}&{A}_{22} \srix 
\smat x_1+x_r \\ x_2\srix +2 B_2 u\right),\\
\mathcal T_8:&=x_2^T \Sigma_2^{-1}k^2x_2 +\sum_{i=1}^m\left(\smat {N}_{i, 21}&{N}_{i, 22} \srix \smat x_1+x_r \\ 
x_2\srix\right)^T \Sigma_2^{-1} \smat {N}_{i, 21}&{N}_{i, 22} \srix \smat x_1+x_r \\ x_2\srix.
\end{align*}}
Exploiting the following equation \begin{align*}
     \frac{d(x_1(t)+x_r(t))}{dt}=&\mat{cc}A_{11} &A_{12}\rix\mat{c}x_1(t)+x_r(t)\\x_2(t)\rix + 2B_1 u(t) \\&+ \sum_{i=1}^m \mat{cc}N_{i,11} 
&N_{i, 12}\rix\mat{c}x_1(t)+x_r(t)\\x_2(t)\rix u_i(t),        
                \end{align*}
we derive  \begin{align*}
   &  \frac{d}{dt}\left((x_1+x_r)^T\Sigma_1^{-1} (x_1+x_r)\right)=2(x_1+x_r)^T\Sigma_1^{-1} 
\frac{d}{dt}(x_1+x_r)\\&=2(x_1+x_r)^T\Sigma_1^{-1}\left(\smat A_{11} &A_{12}\srix \smat x_1+x_r\\x_2\srix + 2B_1 u + \sum_{i=1}^m 
\smat N_{i,11} &N_{i, 12}\srix\smat x_1+x_r\\x_2\srix u_i\right). 
                \end{align*}
Analogously to the above computations the last term can be bounded as follows \begin{align*}
     &\sum_{i=1}^m 2 (x_1+x_r)^T\Sigma_1^{-1} \smat N_{i, 11}&N_{i, 12} \srix \smat x_1+x_r\\ x_2 \srix u_i\\&=\sum_{i=1}^m 2 \left\langle 
\Sigma_1^{-\frac{1}{2}}(x_1+x_r) u_i, \Sigma_1^{-\frac{1}{2}} \smat N_{i, 11}&N_{i, 12} \srix \smat x_1+x_r\\ x_2 \srix \right\rangle_2\\&\leq 
(x_1+x_r)^T\Sigma_1^{-1}(x_1+x_r) k^2+ \sum_{i=1}^m \left(\smat N_{i, 11}&N_{i, 12} \srix \smat x_1+x_r\\ x_2 \srix \right)^T \Sigma_1^{-1}\smat 
N_{i, 11}&N_{i, 12} \srix \smat x_1+x_r\\ x_2 \srix.
                \end{align*}                
This estimate yields  \begin{align*}
     \frac{d}{dt}\left((x_1(t)+x_r(t))^T\Sigma_1^{-1} (x_1(t)+x_r(t))\right)\leq \mathcal T_5+\mathcal T_6.
                \end{align*}
With equation (\ref{derix2}) and the previous steps, we know that \begin{align*}
\frac{d(x_2^T\Sigma_2^{-1} x_2)}{dt}\leq & 2x_2^T\Sigma_2^{-1} \left(\smat A_{21}& A_{22}\srix 
\smat x_1 \\ x_2\srix + B_2u \right)+ x_2^T\Sigma_2^{-1} x_2 k^2\\&+ \sum_{i=1}^m \left(\smat N_{i, 21}&N_{i, 22} \srix \smat x_1\\ x_2 \srix 
\right)^T \Sigma_2^{-1}\smat N_{i, 21}&N_{i, 22} \srix \smat x_1\\ x_2 \srix,
                \end{align*}
such that\begin{align*}
     \frac{d}{dt}\left(x_2(t)^T\Sigma_2^{-1} x_2(t)\right)\leq & \left[\mathcal T_7-2x_2^T\Sigma_2^{-1} ({A}_{21}x_r+B_2 u)\right]\\
&+\left[\mathcal T_8-2 \sum_{i=1}^m (N_{i, 21} x_r)^T \Sigma_2^{-1} \mat{cc} N_{i, 21} & N_{i, 22}\rix \mat{c} x_1 \\x_2\rix \right.\\&\ \ \ \ \ \ 
-\left. \sum_{i=1}^m(N_{i, 21} x_r)^T \Sigma_2^{-1} N_{i, 21} x_r\right].
                \end{align*}           
In summary, we obtain \begin{align} \label{key2}
4 \left\|u(t)\right\|^2_2\geq &\mathcal T_5+\mathcal T_6+\mathcal T_7+\mathcal T_8\\ \nonumber
\geq &\frac{d}{dt} \left((x_1(t)+x_r(t))^T\Sigma_1^{-1} (x_1(t)+x_r(t))\right)+\frac{d}{dt}\left(x_2(t)^T\Sigma_2^{-1} x_2(t)\right)\\ \nonumber
    &+2x_2^T\Sigma_2^{-1} ({A}_{21}x_r+B_2 u)+ 2 \sum_{i=1}^m (N_{i, 21} x_r)^T \Sigma_2^{-1} \smat N_{i, 21} & N_{i, 22}\srix \smat x_1 
\\x_2\srix.   
    \end{align}
Integration of both sides yields  \begin{align}\label{followkey2}
4\int_0^T\left\|u(t)\right\|^2_2 dt\geq & 2\left[\int_0^T x_2^T\Sigma_2^{-1} ({A}_{21}x_r+B_2 u) dt\right.\\&\left.\ \ \ + \int_0^T \sum_{i=1}^m 
(N_{i, 21} x_r)^T \Sigma_2^{-1} \smat N_{i, 21} & N_{i, 22}\srix \smat x_1 
\\x_2\srix dt\right].
                \end{align}
Combining this inequality with (\ref{firstestimate}) leads to \begin{align}\label{errorremovefirst}
\left(\int_0^T \left\|y(t)-y_r(t)\right\|^2_2 dt\right)^{\frac{1}{2}} \leq 2 \sigma \left(\int_0^T\left\|u(t)\right\|^2_2 dt\right)^{\frac{1}{2}}.
\end{align}
 
 The proof for general $\Sigma_2$ relies on the common idea of removing the Hankel singular values step by 
step. The error between the outputs $y$ and $y_r$ can be bounded as follows:\begin{align*}
   \left\|y-y_r\right\|_{L^2_T}\leq   
\left\|y-y_{r_\nu}\right\|_{L^2_T}+\left\|y_{r_\nu}-y_{r_{\nu-1}}\right\|_{L^2_T}+\ldots+\left\|y_{r_2}-y_{r}\right\|_{L^2_T},        
  \end{align*}
where the dimensions $r_i$ of the corresponding states are defined by $r_{i+1}=r_{i}+m(\tilde\sigma_{i})$ for $i=1, 2 \ldots, \nu-1$. Here, 
$m(\tilde\sigma_{i})$ denotes the multiplicity of $\tilde\sigma_{i}$ and $r_1=r$. In the reduction step from $y$ to $y_{r_\nu}$ only the smallest 
Hankel singular value $\tilde\sigma_\nu$ is removed from the system. Hence, by inequality (\ref{errorremovefirst}), we have \begin{align*}            
                             \left\|y-y_{r_\nu}\right\|_{L^2_T}\leq 2 \tilde\sigma_\nu \left\|u\right\|_{L^2_T}.
\end{align*}
Inequality (\ref{errorremovefirst}) can be established as well when comparing the reduced order outputs $y_{r_\nu}$ and $y_{r_{\nu-1}}$. Again, only 
one Hankel singular value, namely $\tilde\sigma_{r_{\nu-1}}$, is removed. At the same time, we have the same kind of inequalities in the ROM as 
before, that are, \begin{align*}
 A_{11}^T \Sigma_1^{-1}+\Sigma_1^{-1}A_{11}+\sum_{i=1}^m N_{i, 11}^T \Sigma_1^{-1} N_{i, 11}+\Sigma_1^{-1} k^2 &\leq -\Sigma_1^{-1}B_1B_1^T
\Sigma_1^{-1},\\
 A_{11}^T \Sigma_1+\Sigma_1 A_{11}+\sum_{i=1}^m N_{i, 11}^T \Sigma_1 N_{i, 11}+\Sigma_1 k^2 &\leq -C_1^TC_1.
                                       \end{align*}
So, by repeatedly applying the above arguments, we obtain\begin{align*}
   \left\|y_{r_j}-y_{r_{j-1}}\right\|_{L^2_T}\leq 2\tilde\sigma_{r_{j-1}} \left\|u\right\|_{L^2_T}
  \end{align*}
for $j=2, \ldots, \nu $. This provides the claimed result.
  \end{proof}
    \end{thm}
Since the bound in Theorem \ref{thmfirsthsv} involves only the sum of distinct diagonal entries of $\Sigma_2$, the result is also true 
when using the sum of all diagonal entries instead. 
\begin{kor}\label{cor2trunhsv}
Let $x(0)=0$, $x_{r}(0)=0$ and $\left\|u(t)\right\|_2\leq k$, $t\in[0, T]$, $T>0$, where $k$ is the constant that enters in 
(\ref{bilinearreachgramtyp2}) and (\ref{bilinearobservgramtyp2}). Then, \begin{align*}
   \left\|y-y_r\right\|_{L^2_T}\leq 2 (\sigma_{r+1}+\sigma_{r+2}+\ldots + \sigma_n)  \left\|u\right\|_{L^2_T},            
  \end{align*}
where $y$ is the output of the original system (\ref{stochstatenewextend}), $y_r$ is the output of the type II BT ROM and 
$\sigma_{r+1}, \ldots,\sigma_n$ are the diagonal entries of $\Sigma_2$.
\end{kor}
The $\mathcal H_\infty$-error of using type II BT depends on the $n-r$ smallest HSVs of the original system. If now only states are neglected that 
correspond  to small values $\sigma_{r+1}, \ldots, \sigma_n$ (hardly reachable and observable states), BT leads to a good approximation of the full 
order system by Corollary \ref{cor2trunhsv}.
\begin{remark}
Theorem \ref{thmfirsthsv} and Corollary \ref{cor2trunhsv} can also be achieved when balancing based on $Q_1$ and $P_2$ satisfying 
(\ref{gengenlyapobs}) and (\ref{newgram2}), respectively. In this case, the same techniques can be used. The key points are inequalities (\ref{key1}) 
and (\ref{key2}). Supposing that $\Sigma_S=\diag(\Sigma_{S,1}, \Sigma_{S, 2})$ is the matrix of HSVs from $Q_1$ and $P_2$, 
(\ref{key1}) then is\begin{align*}
-\left\|y-y_r\right\|^2_2\geq & \frac{d}{dt} (x_1-x_r)^T\Sigma_{S,1} (x_1-x_r)+\frac{d}{dt} x_2^T\Sigma_{S,2} x_2\\&-\left((x_1-x_r)^T\Sigma_{S,1} 
(x_1-x_r)\left\|u\right\|_2^2+x_2^T\Sigma_{S,2} x_2\left\|u\right\|_2^2\right)\\& -2x_2^T\Sigma_{S,2} ({A}_{21}x_r+B_2 u)-2 \sum_{i=1}^m (N_{i, 21} 
x_r)^T \Sigma_{S,2} \smat N_{i, 21} & N_{i, 22}\srix \smat x_1 \\x_2\srix .                                                                           
                  \end{align*}
In order to find (\ref{firstestimate}) from this, we have to guarantee that \begin{align}\label{positive1}
\smat x_1(T)-x_r(T)\\x_2(T)\srix^T\Sigma_{S} \smat x_1(T)-x_r(T)\\x_2(T)\srix \geq \int_0^T \smat x_1-x_r\\x_2\srix^T 
\Sigma_{S}\smat x_1-x_r\\x_2\srix^T\left\|u\right\|_2^2 ds.
                                                                            \end{align}
With the alternative Gramians (\ref{key2}) reads
\begin{align*} 
4 \left\|u(t)\right\|^2_2\geq &\frac{d}{dt} (x_1+x_r)^T\Sigma^{-1}_{S,1} (x_1+x_r)+\frac{d}{dt} x_2^T\Sigma^{-1}_{S,2} 
x_2\\&-\left((x_1+x_r)^T\Sigma^{-1}_{S,1} 
(x_1+x_r)\left\|u\right\|_2^2+x_2^T\Sigma^{-1}_{S,2} x_2\left\|u\right\|_2^2\right)\\& +2x_2^T\Sigma^{-1}_{S,2} ({A}_{21}x_r+B_2 u)+2 \sum_{i=1}^m 
(N_{i, 21} x_r)^T \Sigma^{-1}_{S,2} \smat N_{i, 21} & N_{i, 22}\srix \smat x_1 \\x_2\srix .
    \end{align*}
To obtain (\ref{followkey2}), it needs to be ensured that \begin{align}\label{positive2}
\smat x_1(T)+x_r(T)\\x_2(T)\srix^T\Sigma^{-1}_{S} \smat x_1(T)+x_r(T)\\x_2(T)\srix \geq \int_0^T \smat x_1+x_r\\x_2\srix^T 
\Sigma_{S}^{-1}\smat x_1+x_r\\x_2\srix^T\left\|u\right\|_2^2 ds.
    \end{align}
Both (\ref{positive1}) and (\ref{positive2}) are true if the control $u$ is small enough. Consequently, Theorem \ref{thmfirsthsv} and Corollary 
\ref{cor2trunhsv} hold when using $Q_1$ and $P_2$ instead if $\left\|u(t)\right\|_2\leq k_2$, $t\in[0, T]$, for a sufficiently small constant 
$k_2>0$.
\end{remark}

\section{Conclusions}

We have summerised previous work on balanced truncation for bilinear control systems. We have discussed Gramians that have been studied before and 
how they can 
be used to bound controllability and observability energy functionals, however, these bounds only hold in a small neighbourhood of zero. We proposed 
new Gramians (type II) which contain additional information about the control. With these type II Gramians global energy bounds can be 
found if the controls are assumed to be bounded in a certain way. These bounds justify the use of the alternative Gramians in the context of 
balancing. Based on this motivation, we explained the balancing procedure for bilinear systems which is similar to the one in the 
linear case. Another advantage of using the type II Gramians is the availability of an $\mathcal H_\infty$-error bound for balanced 
truncation of bilinear systems which we proved in this paper. The error bound requires the assumption of having a bounded input to the system. 
Hence, we have overcome the problem of previous works, where no error bound has been established. 

\section*{Acknowledgements}
The author thanks the organisers of the LMS-EPSRC Durham Symposium on Model Order Reduction during which most of the work on this paper was carried 
out.

\bibliographystyle{abbrv}

\end{document}